\documentclass{amsart} 
\usepackage{amssymb,latexsym}
\def\C{\mathbb C}
\def\R{\mathbb R}
\def\Z{\mathbb Z}

\def\N{\mathbb N}
\newtheorem{thm}{Theorem}[section]
\newtheorem{lem}{Lemma}[section]

\newtheorem{prop}{Proposition}[section]
\newtheorem{conj}{Conjecture}[section]

\begin{document}

\title{Equilibrium points of logarithmic potentials on convex domains}
\author{J.K. Langley}
\address{School of Mathematical Sciences
\\University of Nottingham
\\NG7 2RD, UK}
\email{jkl@maths.nott.ac.uk}

\begin{abstract}
Let $D$ be a convex domain in $\C$.
Let $a_k > 0$ be summable constants and let $z_k \in D$. 
If the $z_k$ converge
sufficiently rapidly to $\eta \in \partial D$ 
from within an appropriate Stolz angle
then the function $\sum_{k=1}^\infty a_k /( z  - z_k ) $
has infinitely many zeros in $D$. An example shows
that the hypotheses on the $z_k$ are not redundant, and that two 
recently advanced conjectures are false.\\
M.S.C. 2000 classification: 30D35, 31A05, 31B05.
\end{abstract}
\maketitle

\section{Introduction}

A number of recent papers \cite{CER,ELR,LaRo,LaRo2} have concerned
zeros of functions
\begin{equation}
f(z) = \sum_{k=1}^\infty \frac{a_k}{z - z_k} , 
\label{1}
\end{equation}
and in particular the following conjecture \cite{CER}. 

\begin{conj}[\cite{CER}]\label{conj1}
Let $f$ be given by $(\ref{1})$, where $a_k > 0$ and
\begin{equation}
z_k  \in \C , 
\quad \lim_{k \to \infty } z_k = \infty , \quad \sum_{z_k \neq 0} 
\left| \frac{a_k}{z_k} \right| < \infty .
\label{2}
\end{equation}
Then $f$ has infinitely many zeros in $\C$.
\end{conj}

The assumptions of Conjecture \ref{conj1}
imply that $f$ is meromorphic in the plane
and, assuming that all $z_k$ are non-zero, $f(z)$ is the complex conjugate
of the gradient of the associated subharmonic potential
$u(z) = \sum_{k=1}^\infty a_k \log | 1 - z/z_k | $.
Moreover,
Conjecture \ref{conj1} has a physical interpretation in terms of the 
existence of equilibrium points of the electrostatic field arising
from a system of infinite wires,
each carrying a charge density $a_k$ and perpendicular
to the complex plane at $z_k$ \cite[p.10]{Kellogg}.
Conjecture \ref{conj1} is known to be true
when $\sum_{|z_k| \leq r} a_k = o( \sqrt{r} )$
as $r \to \infty $ \cite[Theorem 2.10]{CER}
(see also \cite[p.327]{GO}), and when
$\inf \{ a_k \} > 0$ \cite{ELR} (see also \cite{LaRo}). 

An analogue of Conjecture \ref{conj1} for a disc
was advanced 
in \cite[Conjecture 2]{borcea}.

\begin{conj}[\cite{borcea}]\label{conj2}
Let $0 < \rho < \infty $ and $\theta \in \R$. 
Let $f$ be given by $(\ref{1})$, where
\begin{equation}
z_k  \in \C ,  \quad |z_k| < \rho, \quad  
\quad \lim_{k \to \infty } z_k = \rho e^{i \theta}, 
\quad a_k > 0, \quad \sum_{k=1}^\infty 
a_k  < \infty .
\label{4}
\end{equation}
Then $f$ has infinitely many zeros in $|z| < \rho $.
\end{conj}
If $f$ satisfies the assumptions of Conjecture \ref{conj2} 
then $\overline{f} = \nabla u$ in $|z| < \rho$, where 
$u(z) = \sum_{k=1}^\infty a_k \log | z - z_k |$. 
Obviously there is no loss of generality in assuming that $\rho = 1$ 
and $\theta = 0$ 
in Conjecture \ref{conj2}. Writing
\begin{equation}
w = \frac1{1-z} , \quad  w_k = \frac1{1 - z_k} , \quad f(z) = w F(w) , 
\label{5}
\end{equation}
where
\begin{equation}
\quad F(w) = \sum_{k=1}^\infty \frac{a_k w_k}{ w - w_k} ,
\label{6}
\end{equation}
it is easy to verify that Conjecture \ref{conj2} is equivalent to
the following.
\begin{conj}\label{conj3}
Let $F$ be given by $(\ref{6})$, where
\begin{equation}
w_k  \in \C ,  \quad {\rm Re \,} w_k > \frac12 ,  \quad  
\quad \lim_{k \to \infty } w_k = \infty , \quad a_k > 0,
\quad \sum_{k=1}^\infty 
a_k  < \infty .
\label{7}
\end{equation}
Then $F$ has infinitely many zeros in ${\rm Re \,} w > 1/2$.
\end{conj}
With the assumptions (\ref{7}), the function $F(w)$ in (\ref{6}) 
is evidently meromorphic in the plane. In \S\ref{example} an 
example satisfying (\ref{6}) and (\ref{7}) will be constructed,
such that $F(w)$ has no zeros in $\C$. Thus Conjectures \ref{conj2} 
and \ref{conj3} are
false, and there is no direct analogue of Conjecture \ref{conj1}
for the unit disc. 

On the other hand the following theorem shows 
in particular that 
if the $z_k$ converge to $\rho e^{i \theta} $ sufficiently rapidly,
and if all but finitely many $z_k$ lie in a sufficiently small Stolz angle,
then the conclusion of Conjecture \ref{conj2} does hold.
It is convenient to state and prove the result when the $z_k$ lie in a
convex domain $D$ and the boundary point 
$\rho e^{i \theta} $ is $1$. There then exists
(see \S\ref{pfthm1}) an open half-plane $H$ such that $D \subseteq H$ and
$1$ lies on the boundary 
$\partial H$, and there is no loss of generality in assuming 
that $H$ is the half-plane ${\rm Re \, } z < 1$. 

\begin{thm}\label{thm1}
Let $D \subseteq \{ z \in \C : {\rm Re \, } z < 1 \}$ be a convex
domain such that $1 \in \partial D$. Let $f$ be given by $(\ref{1})$, where
\begin{equation}
z_k  \in D , \quad a_k > 0, \quad \sum_{k=1}^\infty 
a_k  < \infty .
\label{4'}
\end{equation}
Assume that $1 $ is a
limit point of the set $\{ z_k : k \in \N \}$, and
that there exist real numbers
$\varepsilon > 0$ and $ \lambda  \geq 0$ such that 
\begin{equation}
\sum_{| 1 - z_k| \leq \varepsilon} | 1 - z_k|^\tau < \infty 
\quad \hbox{for all} \quad \tau >  \lambda ,
\label{4''}
\end{equation}
and 
\begin{equation}
\sup \{ | \arg ( 1 - z_k ) | : k \in \N , \, | 1 - z_k| \leq \varepsilon \}
< C ( \lambda   ) = 
\frac{\pi }{2 \lambda }  .
\label{4'''}
\end{equation}
Then there exists a sequence $(\eta_j)$ of zeros of
$f$ satisfying $\eta_j \in D, \lim_{j \to \infty} \eta_j = 1$.
\end{thm}
Note that (\ref{4''})
implies that $\{ z_k : k \in \N \}$ has no limit points $z$ in the 
punctured disc $A$ given by 
$0 < | 1 - z | < \varepsilon $, and that $f$ is meromorphic on $A$.
Moreover, (\ref{4'''}) is obviously satisfied if $\lambda < 1$.

\section{A counterexample to Conjecture \ref{conj3}}\label{example}

Let
\begin{equation}
g(w) = \frac1{w(w-2) ( e^{w-1} + 1 ) } .
\label{8}
\end{equation}
Then $g$ has no zeros, but has simple poles at $0$, $2$ and
\begin{equation}
u_k = 1 + (2k+1) \pi i , \quad k \in \Z .
\label{9}
\end{equation}
Straightforward computations give
\begin{equation}
{\rm Res \, } (g, 0) = \frac{-1}{2(e^{-1} + 1)} = -a , \quad
{\rm Res \, } (g, 2) = \frac{1}{2(e + 1)} = b , 
\label{10}
\end{equation}
and, using (\ref{9}), 
\begin{equation}
{\rm Res \, } (g, u_k) = \frac{-1}{u_k(u_k-2)} = 
\frac{-1}{(u_k - 1)^2 - 1}  
= \frac1{ (2k+1)^2 \pi^2 + 1} = c_k .
\label{12}
\end{equation} 
Then $b$ and the $c_k$ evidently satisfy
\begin{equation}
b > 0, \quad c_k > 0, \quad \sum_{k \in \Z} c_k < \infty .
\label{13}
\end{equation}
Next, let
\begin{equation}
h(w) = - \frac{a}{w} + \frac{b}{w - 2}  + 
\sum_{k \in \Z} \frac{c_k}{w - u_k} , \quad
L(w) = h(w) - g(w) .
\label{14}
\end{equation}
By (\ref{8}), (\ref{9}), (\ref{10}), (\ref{12}) and
(\ref{13}) the function $h(w)$ is meromorphic in the plane, 
and $L(w)$ is an entire function.

Let $m$ be a large positive integer, let $R = 4 m \pi $, and use $c$
to denote positive constants independent of $m$. Then
simple estimates give
\begin{equation}
|g(w)| \leq \frac{c}{R^2} 
\quad \hbox{for} \quad |w - 1| = R
\label{15}
\end{equation} 
and, as $m \to \infty$, 
\begin{equation}
|h(w)| \leq  \frac{c}{R}  
+ c \sum_{k \in \Z , |k| \geq m } c_k 
+ c \sum_{k \in \Z , |k| <  m } \frac{c_k}{R} = o(1)  
\quad \hbox{for} \quad |w - 1 | = R .
\label{16}
\end{equation} 
Combining (\ref{15}) and (\ref{16}) shows that $L(w) \equiv 0$ in
(\ref{14}), so that $h = g$ has no zeros, and applying 
the residue theorem in conjunction with (\ref{15}) now gives
\begin{equation}
a =  b + \sum_{k \in \Z} c_k .
\label{17}
\end{equation}
Hence $h(w)$ may be expressed using (\ref{17}) in the form 
\begin{eqnarray}
h(w) &=& 
b \left( \frac1{w-2} - \frac1w \right) + 
\sum_{k \in \Z} c_k \left( \frac1{w - u_k} - \frac1w \right) 
\nonumber \\
&=& \frac1w \left( \frac{2b}{w - 2} + 
\sum_{k \in \Z} \frac{c_k u_k }{w - u_k} \right) .
\label{18}
\end{eqnarray}
By (\ref{9}), (\ref{13}) and (\ref{18}) the function $F(w) =
wh(w)$ may be written in the
form 
\begin{equation}
F(w)  = \sum_{k=1}^\infty \frac{d_k v_k}{w - v_k},
\quad {\rm Re \, } v_k \geq 1, \quad 
v_k \to \infty, \quad d_k > 0, \quad \sum_{k=1}^\infty d_k < \infty .
\label{18'}
\end{equation}
Here $F$ evidently satisfies the requirements of
(\ref{6}) and (\ref{7}), but $F$ has no zeros in $\C$,
since $h$ has no zeros and $h(0) = \infty$.
\\\\
{\em Remark.} It is conjectured further in \cite[Conjecture 6]
{borcea} that if $f$ satisfies
(\ref{1}) and (\ref{2}) with $a_k \overline{z}_k > 0$ for each $k$ 
then $f$ has infinitely many zeros in $\C$. The example (\ref{18'}), 
with $a_k = d_k v_k$ and $a_k \overline{v}_k = d_k |v_k|^2 > 0$,
shows that this conjecture is also false.

\section{An auxiliary result needed for Theorem \ref{thm1}}

The proof of Theorem \ref{thm1} rests upon the following proposition,
which  
concerns
functions in the plane of the form (\ref{6}), and uses standard notation
from \cite[p.42]{Hay2}. 

\begin{prop}\label{prop1}
Let $0 < \sigma \leq 1$.
Let $F$ be given by 
$(\ref{6})$, where
\begin{equation}
w_k  \in \C ,  \quad {\rm Re \,} w_k > 0 ,  \quad  
\quad a_k > 0,
\quad \sum_{k=1}^\infty 
a_k  < \infty .
\label{7'}
\end{equation}
Assume that the set $\{ w_k : k \in \N \} $ is unbounded and that
there exist real numbers $R > 0$ and
$\lambda \geq 0$ such that 
\begin{equation}
\sum_{|w_k| \geq R} |w_k|^{- \tau } < \infty 
\quad \hbox{for all} \quad \tau > \lambda ,
\label{7''}
\end{equation}
and 
\begin{equation}
s = \sup \{ | \arg w_k | : k \in \N , \, |w_k| \geq R \} <  
C ( \lambda , \sigma ) = 
\frac{2}{\lambda } \arcsin \sqrt{ \frac{\sigma}{2} } .
\label{t1}
\end{equation}
Then there exists a transcendental
meromorphic function $G$ with
\begin{equation}
F(w) = G(w)(1 + o(1)) \quad \hbox{as} \quad w \to \infty ,
\label{t0}
\end{equation}
and the Nevanlinna deficiency
$\delta (0, G)$ of the zeros of $G$ satisfies 
$\delta (0, G) < \sigma $. In particular, $F(w)$ has a sequence of 
zeros tending to infinity.
\end{prop}

The zero-free example of (\ref{18'}) has $\lambda = 1$ and
$\delta (0, F) = \sigma = 1$, and all  
its poles lie in ${\rm Re \, } w \geq 1$, so that
Proposition \ref{prop1} is essentially sharp. 

To prove Proposition \ref{prop1}, 
assume that $F$ is as in the statement of Proposition \ref{prop1}.
It follows from (\ref{7''}) that the set $\{ w_k : k \in \N \}$
has no limit points $w$ with $R < |w| < \infty $. In particular,  
$F$ is meromorphic in the region $2R \leq |w| <
\infty$ with an essential singularity at infinity. The existence
of a transcendental meromorphic function $G$ satisfying (\ref{t0})
then follows from a result of Valiron \cite[p.15]{Valiron}
(see also \cite[p.89]{Bie}). In particular, $G$ is constructed 
\cite{Valiron} so that
$F$ and $G$ have the same
poles and zeros in $|w| \geq 2R$. If $|w| \geq 4R$ then (\ref{7'}) gives
$$
|F(w)| \leq |F_1(w)| + O(1), \quad 
F_1(w) =  \sum_{|w_k| \geq 2R} \frac{a_k w_k}{ w - w_k} ,
$$
so that
$$
m(r, G) \leq m(r, F_1) + O(1) = O(1) 
$$
as $r \to \infty$, by \cite[p.327]{GO}. 
Since the poles $w_k$ of $G$ have exponent of
convergence at most $\lambda$ by (\ref{7''}), it follows that $G$ has
lower order $\mu \leq \lambda$. 

Choose $s_0, s_1, s_2$ with
\begin{equation}
s < s_0  < s_1 < s_2 <  \min \{ \pi ,  C( \lambda , \sigma ) \} ,
\label{sj}
\end{equation} 
where $s$ is as in (\ref{t1}) and satisfies $s \leq \pi /2 $ by
(\ref{7'}).
The proof of Proposition \ref{prop1}
requires the following two lemmas.

\begin{lem}\label{lem0}
The function $F$ satisfies 
$\liminf_{r \in \R, r \to + \infty } r |F(-r)| > 0 $.
\end{lem}

\begin{proof}
Let $r > 0$ and write $w_k = u_k + i v_k$ with $u_k$ and 
$v_k$ real. 
Let
$$
p_k (r) = {\rm Re \, } \left( \frac{w_k}{r + w_k} \right) =
\frac{ u_k (r + u_k ) + v_k^2 }{ (r + u_k)^2 + v_k^2 } .
$$
Then (\ref{7'}) gives $p_k (r) > 0$ and there exists $d > 0$ such that
$p_1(r) > d/r$ as $r \to \infty$. Hence, again as $r \to \infty$,
$$
r |F(-r)|  \geq - r {\rm Re \,}  F(-r)  
= r \sum_{k=1}^\infty a_k p_k (r) \geq r a_1 p_1 (r) > a_1 d .
$$
\end{proof}
\begin{lem}\label{lem1}
There exists $M_1 > 0$ such that $|F(w)| \leq M_1$
for all large $w$ lying outside the region 
$| \arg w | < s_0 $. 
\end{lem}

\begin{proof} 
This follows from (\ref{7'}), (\ref{t1}) and (\ref{sj}), since there exists
a positive constant $M_2 $ such that 
$|w - w_k| \geq M_2 |w_k| $ for all such $w$ and all $k \in \N$.
\end{proof}

The proof of Proposition \ref{prop1} may now be completed using
Lemmas \ref{lem0} and \ref{lem1}. 
Assume that 
$\delta (0, G) \geq \sigma$. Then
Baernstein's
spread theorem \cite{baernstein} gives
a sequence $r_m \to \infty$
and, for each $m$, a subset $I_m$ of the circle $|w| = r_m $, of 
angular measure at least 
$$
\min \left\{ 2 \pi ,
\frac{4}{\mu} \arcsin \sqrt{ \frac{\sigma}{2} }  \right\} - o(1) \geq
\min \{ 2 \pi , 2 C( \lambda , \sigma  ) \}   - o(1) \geq 2 s_2 ,
$$
using (\ref{t1}) and (\ref{sj}), and such that
\begin{equation}
\lim_{m \to \infty } 
\frac{ \max \{ \log | G(w) | : \quad w  \in I_m \} }{ \log r_m } = - \infty. 
\label{t2}
\end{equation}
Let $m$ be large, and consider the 
function $v(w) = \log |F(w)|$, which is subharmonic on the domain 
$$
\Omega = \{ w \in \C : r_m /4 < |w| < r_m, 
s_0 < \arg w < 2 \pi - s_0 \} .
$$
Then $v$ is bounded above on $\Omega$, by Lemma
\ref{lem1}. But the intersection $J_m$ of $I_m$ with
the arc $\{ w \in \C : |w| = r_m, s_1 < \arg z < 2 \pi - s_1 \}$ has
angular measure at least $2( s_2 - s_1 )$, so that
standard estimates for the harmonic measure
of $J_m$ at $- r_m /2$ now give
\begin{equation}
\omega (-r_m /2, J_m, \Omega ) \geq M_3 > 0,
\label{t3}
\end{equation}
where $M_3$ is independent of $m$. Since (\ref{t0}) implies that
(\ref{t2}) holds with $G$ replaced by $F$,
combining Lemma \ref{lem1} 
with (\ref{t3})
and the two-constants theorem \cite[p.42]{Nev} leads to
\begin{equation}
r_m F( - r_m /2 )  \to 0 \quad \hbox{as}  \quad m \to \infty .
\label{t4}
\end{equation}
But (\ref{t4}) contradicts Lemma \ref{lem0}, and this completes the
proof of Proposition \ref{prop1}.

\section{Proof of Theorem \ref{thm1}}\label{pfthm1}
Assume that $f$ and $D$
satisfy the hypotheses of Theorem \ref{thm1}. 
Define $F$ using the transformations (\ref{5}) and (\ref{6}). 
Then $F$ satisfies
the hypotheses of Proposition \ref{prop1} with $R = 1/\varepsilon$
and $\sigma = 1$. 
Thus $F$ has a sequence of zeros tending to infinity, and so $f$
has a sequence $(\eta_j)$ of zeros with 
$\lim_{j \to \infty} \eta_j = 1$.

It remains only to show that such a sequence $(\eta_j)$ exists with, in
addition, $\eta_j \in D$, and this is done by a 
standard argument of Gauss-Lucas type. 
Let $\eta = \eta_j$ with $j$ large, and assume that $\eta \not \in D$. 
Since $D$ is convex 
the supremum and infimum of $\arg (z - \eta)$ on $D$ differ by at most $\pi$.
Hence there exist an open half-plane $H$, with $D \subseteq H$ and
$\eta \in \partial H$, and a linear transformation $u = T(z) =  ( z - \eta )/a$
mapping $H$ onto ${\rm Re \, u} > 0$. Writing $u_k = T(z_k)$ then gives
$$
0 = {\rm Re \, } (a f(\eta )) = - {\rm Re \, }
\left( \sum_{k=1}^\infty \frac{a_k}{ u_k} \right) < 0 .
$$
This contradiction completes the proof of Theorem \ref{thm1}.

\end{document}